\documentclass{elsarticle}

\usepackage{fullpage}
\usepackage{amsmath}
\usepackage{amssymb}
\usepackage{amsthm}
\usepackage{graphicx}
\usepackage{color}
\newcommand{\Vol}{\text{Vol}}

\newtheorem{theorem}{Theorem}[section]
\newtheorem{lemma}[theorem]{Lemma}
\newtheorem{defn}[theorem]{Definition}
\newtheorem{corollary}[theorem]{Corollary}

\begin{document}

\begin{abstract}
We answer a problem posed recently by Knuth: an $n$-dimensional box, with edges lying on the positive coordinate axes and generic edge lengths $W_1 < W_2 < \cdots < W_n$, is dissected into $n!$ pieces along the planes $x_i = x_j$.  We describe which pieces have the same volume, and show that there are $C_n$ distinct volumes, where $C_n$ denotes the $n$th Catalan number.

\end{abstract}

\title{A solution to one of Knuth's permutation problems} 
\author{Benjamin Young}
\address{byoung@math.mcgill.ca}
\maketitle

\section{Introduction}

In a recent talk~\cite{knuth-2010}, D. Knuth posed the following problem. Consider the $n$-dimensional box $ B= [0, W_1] \times \cdots \times [0, W_n]$, where $W_1 < W_2 < \cdots < W_n$. If $\pi$ is a permutation in $S_n$, the symmetric group on $n$ letters, define the region
\[
C_{\pi} = \{x \in B |x_{\pi(1)} \geq x_{\pi(2)} \geq \cdots \geq x_{\pi(n)}\}. 
\]
In other words, we dissect $B$ by cutting it along the planes $x_i = x_j$, for $1 \leq i < j \leq n$.  Each $C_{\pi}$ is a piece of this dissection.  Let us view the volume of $C_{\pi}$ as a polynomial in the $W_i$.  How many distinct volumes are there amongst the $C_{\pi}$, and which $C_{\pi}$ have the same volume?  

\begin{figure}
\caption{Knuth's problem in dimension 3}
\label{fig:dim3_version}
\begin{center}
\includegraphics[height=1.3in]{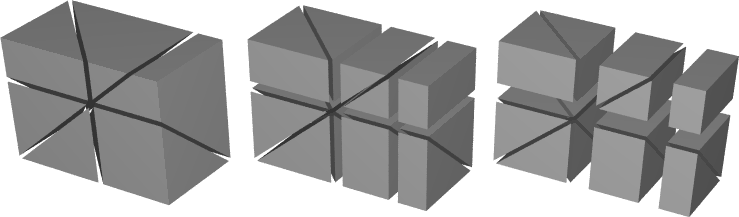}
\end{center}
\end{figure}

See Figure~\ref{fig:dim3_version} for the case $n=3$, in which $C_{132}$ and $C_{123}$ have the same volume and all others have distinct volumes.  The left-hand image shows the original problem; the other two images show $B$ being dissected further along the planes $x_i = W_j$, so that the volumes may be more easily computed. 

\begin{defn}
\label{defn:psi}
Let $\mathcal{P}$ denote the set of all partitions.  Let $\pi$ be a permutation with matrix $[a_{ij}]$.  Define $\psi:S_n \rightarrow \mathcal{P}$ to be the map which sends $\pi$ to the partition whose Young diagram is 
\[\{(i', j')\;:\; a_{ij}=0 \text{ for all }i \leq i'\text{ and }  j \leq j' \}.  \]
\end{defn}
In other words, we cross out all matrix entries which lie weakly below and/or to the right of every one in the permutation matrix for $\pi$ (see Figure~\ref{fig:essentialset}, center image).  The entries which are not crossed out form the Young diagram of $\psi(\pi)$.
\begin{theorem}
\label{thm:main}
If $\pi$ and $\sigma$ are permutations, then
$\Vol(C_{\pi}) = \Vol(C_{\sigma})$ if and only if $\psi(\pi) = \psi(\sigma)$.
\end{theorem}

We defer the proof of this theorem to the end of the paper.  However, there is an immediate corollary, if we appeal to a few results in the literature: 

\begin{corollary}
\label{cor:main}
The number of distinct elements of the set $\{\Vol(C_{\pi}) \;:\;\pi \in S_n\}$ is $C_n = \frac{1}{n+1}\binom{2n}{n}$, the $n$th Catalan number.
\end{corollary}

\begin{proof}
Observe that $\psi(\pi)$ is closely related to a well-known construction, namely that of the \emph{diagram} of the permutation $\pi$.  To construct the diagram of $\pi$, one crosses out all entries \emph{directly} below and \emph{directly} to the right of each of the ones in the matrix for $\pi$.  The result need not be a Young diagram (see Figure~\ref{fig:essentialset}, right image).
 As observed by Reifergerste~\cite{reifergerste-2003}, this procedure yields a Young diagram (and hence coincides with our $\psi(\pi)$) precisely when $\pi$ is 132-avoiding.  In other words, our $\psi$ map yields precisely the rank-zero piece of Fulton's essential set~\cite{fulton-1992, eriksson-linusson-1996}; the entire essential set has rank zero precisely when $\pi$ is 132-avoiding.  
Alternatively, one can see directly that boundary of the Young diagram for $\psi(\pi)$ is always a Dyck path~\cite{fulmek-2003}.  Both 132-avoiding permutations and Dyck paths are enumerated by the Catalan numbers.
\end{proof}

We do not know of a good reason why this problem, or our solution, should have anything to do with combinatorial representation theory; the map $\psi$ as defined above arises naturally in our solution.  

We note that Knuth's original setting of the problem~\cite{knuth-2010} is slightly different.
Namely, fix weights $W_1 < \cdots < W_n$, and let $X_1, \ldots, X_n$ be uniform random variables on $[0,1]$.  We rank the quantities $x_i = W_iX_i$ from smallest to largest.  If $\pi$ is a permutation on $n$ letters, define the event
$E_{\pi}: x_{\pi(1)} \geq x_{\pi(2)} \geq \cdots \geq x_{\pi(n)}.$
Knuth observed that when $n \geq 3$, certain of these events $E_{\pi}$ occur with the same probability regardless of the choice of $W_i$.  Theorem~\ref{thm:main} now classifies the events $E(\pi)$ which occur with the same probability.

We would like to thank D. Knuth for helpful correspondence.

\section{A refinement of the dissection}

We will proceed in the manner suggested in Figure~\ref{fig:dim3_version}: we subdivide the box $B$ further, along the hyperplanes $x_i = W_j$.  Once this is done, all pieces have very simple shapes, and are easily understood.

\begin{defn} Let $W_0=0$, and define 
\begin{align*}
a_i &= W_i - W_{i-1} > 0 \text{ for }1 \leq i \leq n, \\
\mathcal{B} &= \{1,2, \ldots, n\}^n, \\
I &= \{1\} \times \{1,2\} \times \{1,2,3\} \times \cdots \times \{1,2,3,\ldots,n \} \subseteq \mathcal{B}
\end{align*}
and for $\rho = (\rho_1, \ldots, \rho_n) \in \mathcal{B}$, we define the open box
\[
B_{\rho} = (W_{\rho_1-1}, W_{\rho_1}) \times (W_{\rho_2-1}, W_{\rho_2}) \times \cdots \times (W_{\rho_n-1}, W_{\rho_n}).
\]
\end{defn}
Note that the dimensions of $B_{\rho}$ are $a_{\rho_1} \times a_{\rho_2} \times \cdots \times a_{\rho_n}$. Observe that those boxes $B_\rho$ for which $\rho \in I$ lie within $B$, and indeed partition $B$ up to a set of volume zero (namely, the boundaries of the boxes).  
Also, note that if $\rho \in \mathcal{B}$ and $\rho_i = \rho_j$ for some $i \neq j$, then $B_{\rho}$ is symmetric about the hyperplane $\{x_i = x_j\}$,  whereas if $\rho_i < \rho_j$, then the hyperplane $\{x_i = x_j\}$ does not intersect $B_{\rho}$ at all.  

The motivation for all of these definitions is to simplify the computation of the volumes of the $C_{\pi}$.  We begin with the following immediate observation:

\begin{lemma}
\label{lem:trivial_intersections}
For any $\rho \in \mathcal{B}$ and any $\pi \in S_n$, $\rho_{\pi(1)} \geq \rho_{\pi(2)} \geq \cdots \rho_{\pi(n)}$ if and only if all points $x \in B_{\rho}$ satisfy $x_{\pi(1)} \geq x_{\pi(2)} \geq \cdots x_{\pi(n)}.$
\end{lemma}

The symmetric group $S_n$ acts on $\mathcal{B}$ by permuting coordinates.  Each box $B_{\rho}$ has a stabilizer $G_{\rho} \leq S_n$ under this action.  In fact, $G_{\rho}$ is isomorphic to a product of symmetric groups
\[
G_{\rho} \simeq S_{n_1} \times \cdots \times S_{n_k}
\]
where $n_j$ is the number of occurrences of the number $j$  in $\rho$.  Observe that $G_{\rho}$ also acts faithfully on $B_{\rho}$ by permuting coordinates, and so partitions $B_{\rho}$ into $|G_{\rho}|$ equal-volume fundamental domains.  We thus have the following volume computation:
\begin{lemma}
\label{lem:volume_computation}
\[
\Vol(C_{\pi} \cap B_{\rho}) = 
\begin{cases}
0 & \text{ if } C_{\pi} \cap B_{\rho} = \emptyset, \\
\frac{1}{|G_{\rho}|}a_{\rho_1}a_{\rho_2}\cdots a_{\rho_n} &\text{ otherwise.}
\end{cases}
\]
\end{lemma}

\section{Proof of the main theorem}
For the following lemmata and their proofs, we adopt the following notation: 
Let $\rho = (\rho_1, \ldots, \rho_n) \in \mathcal{B}$, $\pi \in S_n$, and let $\lambda = (\lambda_1, \lambda_2, \cdots, \lambda_n) \in \mathcal{B}$ be such that $\lambda_i = \rho_{\pi(i)}$.

\begin{lemma}
\label{lem:partition_meet_condition}
$C_{\pi}$ meets $B_{\rho}$ if and only if $\lambda$ is a partition and $\lambda(i) \geq \pi(i)$.
\end{lemma}

\begin{proof}
By Lemma~\ref{lem:trivial_intersections}, $C_{\pi}$ meets $B_{\rho}$ if and only if $\rho \in I$ and $\rho_{\pi(1)} \geq \cdots \geq \rho_{\pi(n)}$.  Now, $\rho \in I \Leftrightarrow \rho_i \leq i \Leftrightarrow \lambda_i \leq \pi(i)$; similarly, $\rho_{\pi(1)} \geq \cdots \rho_{\pi(n)}$ is equivalent to $\lambda_1 \geq \cdots \geq \lambda_n$.
\end{proof}

Recall that the set of integer partitions forms a distributive lattice, \emph{Young's lattice}, under the partial order of inclusion of Young diagrams. See, for example, \cite[Section 7.2]{stanley-2001} for an introduction to Young's lattice. 

\begin{defn}
Let $\lambda^{\max}(\pi) =  \bigcup \{\mu \in \mathcal{P}\;:\; \mu \text{ is a partition with $n$ parts and } \mu_i \leq \pi(i)\}$, where $\bigcup$ denotes union of Young diagrams (the least upper bound in Young's lattice).
\end{defn}

\begin{lemma}
\label{lem:lambda_max_criterion}
$C_{\pi}$ meets $B_{\rho}$ if and only if $\lambda$ is a partition and $\lambda \subseteq \lambda^{\max}$ as Young diagrams.
\end{lemma}
\begin{proof}
It is easy to check that if $\lambda$ and $\mu$ are partitions which meet the condition of Lemma~\ref{lem:partition_meet_condition}, then so is $\lambda \cup \mu$ (their union as Young diagrams).  Moreover, if $\nu \subseteq \lambda$, then $\nu$ meets the conditions of Lemma~\ref{lem:partition_meet_condition}.  As such, the condition of Lemma~\ref{lem:partition_meet_condition} is equivalent to $\lambda \subseteq \lambda^{\max}$.
\end{proof}

\begin{figure}
\caption{Three of the constructions described in this article, applied to the permutation $\pi=42531$.  From left to right: $\lambda^{\max}(\pi) = (4,2,2,2,1)$, $\psi(\pi) = (3,1,1,1)$ , and the diagram of $\pi$.}
\label{fig:essentialset}
\begin{center}
\includegraphics[height=1.2in]{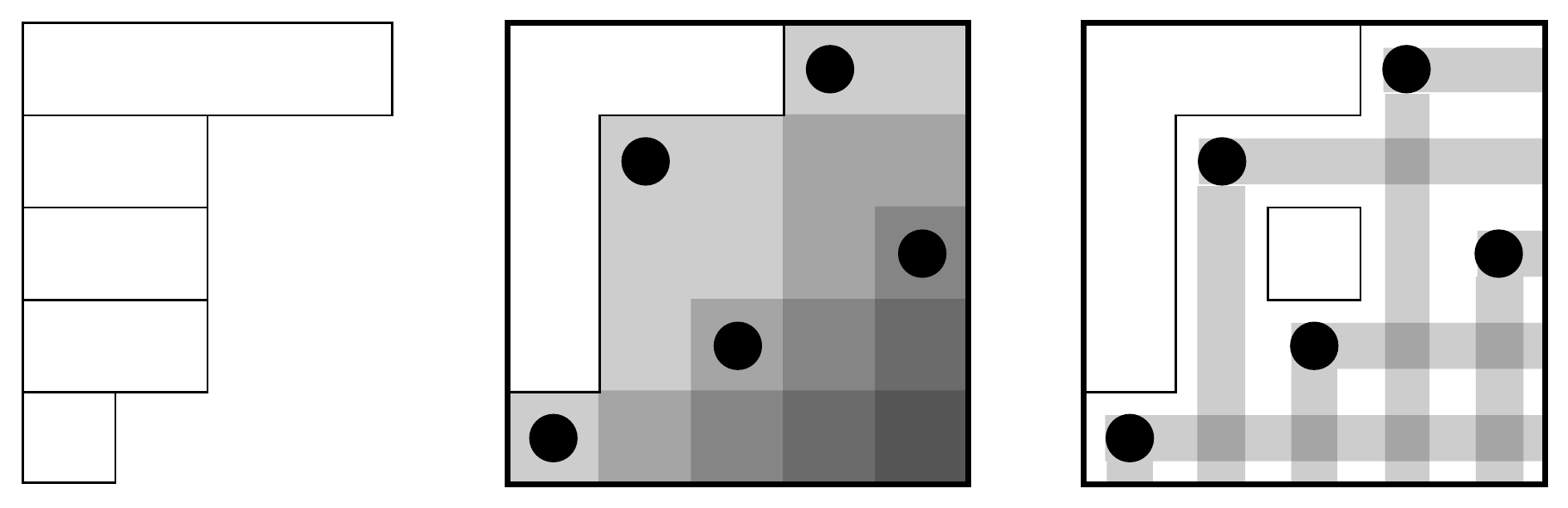}
\end{center}
\end{figure}

\begin{proof}[Proof of Theorem~\ref{thm:main}:]  If $\lambda$ is a partition, write $\rho(\lambda) = (\lambda_{\pi^{-1}(1)}, \ldots, \lambda_{\pi^{-1}(n)}$). Taking $\rho = \rho(\lambda)$ and applying Lemmas~\ref{lem:volume_computation}~and~\ref{lem:lambda_max_criterion}, we see that
\begin{align*}
\Vol(C_{\pi}) 
&= \sum_{\lambda \subseteq \lambda^{\max}(\pi)} \frac{1}{|G_{\rho(\lambda)}|} \prod_i a_{\lambda_i} \\
&= \sum_{\lambda \subseteq \lambda^{\max}(\pi)} \frac{1}{|G_{\lambda}|} \prod_i a_{\lambda_i}.
\end{align*}
The latter equality holds because $G_{\rho}$ is isomorphic to $G_{\sigma \cdot \rho}$ for any permutation $\sigma \in S_n$.  As such, $\Vol(C_{\pi}) = \Vol(C_{\pi'})$ if and only if $\lambda^{\max}(\pi) = \lambda^{\max}(\pi')$.  

Next, we need a concrete description of $\lambda_{\max}(\pi)$.  Let $\lambda$ be a partition such that $\lambda_i \leq \pi(i)$.  In particular, 
\begin{align*}
\lambda_1 &\leq \pi(1), \\
\lambda_2 &\leq \min\{\lambda_1, \pi(2)\} \leq \min\{\pi(1), \pi(2)\}, \\
&\vdots \\
\lambda_n & \leq \min\{\lambda_{n-1}, \pi(n)\} \leq \min\{\pi(1), \ldots, \pi(n)\}.
\end{align*}
Now, $\lambda$ is maximal in Young's lattice if we replace all of the above inequalities with equalities.  Therefore, $\lambda^{\max}_i = \min\{\pi(1), \ldots, \pi(i)\}$.  

Recalling Definition~\ref{defn:psi}, we now compare $\lambda^{\max}(\pi)$ to $\psi(\pi)$. Observe that the permutation matrix for $\pi$ has ones in positions $(i, \pi(i))$ and zeros elsewhere, so 
the $i$th part of $\psi(\pi)$ is $\min\{\pi(1), \pi(2), \ldots, \pi(i)\} - 1$.  In other words, one obtains $\psi(\pi)$ by deleting the first column of the Young diagram of $\lambda^{\max}$; this column is necessarily of height $n$, so one can also reconstruct $\lambda^{\max}(\pi)$ given $\psi(\pi)$ (see Figure~\ref{fig:essentialset}, left and center images).
We conclude that if $\pi, \pi'$ are permutations in $S_n$, then 
\[\Vol(C_{\pi}) = \Vol(C_{\pi'}) 
\Leftrightarrow \lambda^{\max}(\pi) = \lambda^{\max}(\pi')
\Leftrightarrow \psi(\pi) = \psi(\pi').\]
\end{proof}

\bibliographystyle{plain}
\bibliography{knuth}

\end{document}